\documentclass[12pt,a4paper]{article}

\usepackage{amsmath, amssymb}
\usepackage{accents}
\usepackage{float}
\usepackage{hyperref}
\usepackage[linesnumbered,ruled,vlined]{algorithm2e}

\providecommand{\DontPrintSemicolon}{\dontprintsemicolon}
\usepackage{tikz}
\usepackage{multicol}
\usepackage{bm}
\usepackage{xcolor,colortbl}
\usepackage{subcaption}
\usepackage{caption}

\usepackage{multicol}
\usepackage{graphicx}
\usepackage[left=2.5cm, right=2.5cm, top=2.5cm, bottom=2.5cm]{geometry}
\usepackage[affil-it]{authblk}
\usepackage{fancyhdr}
\usepackage{abstract}
\usepackage[backend=bibtex,style=phys]{biblatex}
\usepackage{pifont}

\DeclareMathOperator*{\argmin}{arg\,min}

\title{\textbf{Parallel Trust-Region Approaches in Neural Network Training: Beyond Traditional Methods}}
\author[1]{Ken Trotti}
\author[1,2]{Samuel A. Cruz Alegr\'{i}a}
\author[1,3]{Alena Kopani\v{c}\'{a}kov\'{a}}
\author[1,2]{Rolf Krause}
\affil[1]{Universit\`{a} della Svizzera italiana, Lugano, Switzerland}
\affil[2]{UniDistance Suisse, Brig, Switzerland}
\affil[3]{Brown University, Providence, USA}
\date{\today} 

\begin{document}

\maketitle
\begin{abstract}
\noindent 
We propose to train neural networks (NNs) using a novel variant of the ``Additively Preconditioned Trust-region Strategy'' (APTS). The proposed method is based on a parallelizable additive domain decomposition approach applied to the neural network's parameters.
Built upon the TR framework, the APTS method ensures global convergence towards a minimizer.
Moreover, it eliminates the need for computationally expensive hyper-parameter tuning, as the TR algorithm automatically determines the step size in each iteration.
We demonstrate the capabilities, strengths, and limitations of the proposed APTS training method by performing a series of numerical experiments. The presented numerical study includes a comparison with widely used training methods such as SGD, Adam, LBFGS, and the standard TR method. 
\end{abstract}
%

\section{Introduction}
We consider the following \textit{supervised learning problem}:
\begin{equation}\label{eq:min_problem}
\min_{\theta\in\mathbb{R}^n} \ell(\theta; \mathcal{D}) \;:=\; \min_{\theta\in\mathbb{R}^n} \frac{1}{p} \sum_{i=1}^p l(\mathcal{N}(x_i; \theta, \mathcal{D}), c_i),
\end{equation}
where $\ell \colon \mathbb{R}^n \longrightarrow \mathbb{R}$ is the \textit{training loss}, $\mathcal{N}$ denotes the neural network (NN) and $\theta\in\mathbb{R}^n$ are its parameters. 
The tuple $(x_i,c_i) \in \mathcal{D}$ is an input-target pair in the labeled dataset  $\mathcal{D}$. 
The \textit{loss function} $l \colon \mathbb{R}_{\geq0}^d \times \mathbb{R}_{\geq0}^d \longrightarrow \mathbb{R}$ measures the difference between the predicted target $\hat{c}_i:=\mathcal{N}(x_i; \theta,\mathcal{D})$ and the exact target $c_i\in\mathbb{R}^d$~\cite{goodfellowEtAl2016}. 
Minimizing \eqref{eq:min_problem} is called \textit{training} and successful training implies the model's capability to accurately predict unseen data~\cite{murphy2012}, i.e., the data not contained in the  dataset $\mathcal{D}$.

Several optimization techniques are applicable for training neural networks. 
The most prominent approaches include for example stochastic gradient descent (SGD) \cite{zinkevichEtAl2010, bottouNocedal2018, daiZhu2018}, ADAptive Moment estimation (Adam) \cite{kingmaAdam2014, reddiEtAl2018OpenReview, ruder2016}, and the limited-memory Broyden Fletcher Goldfarb Shanno (LBFGS) algorithm \cite{nocedalWright1999}, each with unique features and challenges. 
The SGD method is efficient, and generalizes well, but requires careful hyper-parameter tuning and diminishing step sizes for global convergence~\cite{zinkevichEtAl2010,bottouNocedal2018,daiZhu2018}. Adam is efficient and adjusts its learning rate adaptively but lacks theoretical global convergence guarantees for non-convex problems and may even face issues with convex problems \cite{kingmaAdam2014,reddiEtAl2018OpenReview}. 
The LFBGS method approximates the Hessian using the secant pairs~\cite{nocedalWright1999} and it is therefore suited for solving ill-conditioned problems. 
However, the LBFGS method is computationally more expensive than the first-order methods and it requires a use of novel techniques for evaluating the secant pairs in stochastic settings~\cite{berahasNocedalTakac2016,berahasEtAl2022}.
Trust-Region (TR) methods~\cite{connGouldToint2000}, though less common, are gaining recognition in machine learning.
Their popularity arises mostly from the fact that they adjust the step size dynamically at each iteration and that they offer global convergence convergences for convex as well as non-convex minimization problems~\cite{nocedalWright1999,connGouldToint2000}. 
However, its extension to multi-level settings \cite{grossKrause2009TR} and to stochastic settings requires non-trivial theoretical and algorithmic adjustments, see for example~\cite{rafatiDeguchyMarcia2018,curtisScheinbergShi2019,gaoNg2022,grattonJeradToint2023, gratton2023multilevel,kopanickovaKrause2022}.



With growing data and network complexity, parallelization approaches become essential~\cite{catlett1991, blum1988, srivastavaEtAl2015, grosse2018}. 
Current parallel approaches include data-parallel and model-parallel methods. 
Data parallelism divides data among processing units, each training a model copy, which is especially useful for large datasets \cite{shallue2018}. 
This paradigm has been applied in SGD variants \cite{zinkevichEtAl2010, nicholsEtAl2021, zhangEtAl2015}, ensemble learning \cite{hansenSalamon1990, dietterich2000}, and federated learning \cite{mcmahanEtAl2017, mishchenkoEtAl2022, malinovskyEtAl2022}. 
In contrast, model parallelism aims at training different parts of the NN across different processing units \cite{benNun2019}.

Recently, domain decomposition (DD) methods show promise in parallelizing NN's training~\cite{shallue2018, benNun2019, nicholsEtAl2021}. 
DD methods were originally proposed for solving discretized partial differential equations. 
The idea behind these methods is to subdivide the main problem into smaller subproblems, a solution of which is used to enhance the solution process of the original problem. 
The subproblems are typically solved independently of each other, i.e., in an additive manner, which in turn enables parallel computing \cite{chanEtAl1994, toselliWidlund2004, gander2006, erhelEtAl2014, mathew2008}. 
Notable examples of additive nonlinear DD methods include ASPIN~\cite{caiKeyes2002}, GASPIN~\cite{grossKrause2021}, RASPEN~\cite{doleanEtAl2016}, or APTS~\cite{grossKrause2009}. 
 While most of these methods were developed for solving nonlinear systems of equations, APTS and GASPIN were specifically designed for addressing non-convex minimization problems, such as one given in Equation~\eqref{eq:min_problem}. 
Furthermore, both APTS and GASPIN employ a TR globalization strategy and, therefore, provide global convergence guarantees.

In the field of machine learning, DD approaches are increasingly utilized to enhance NN training. The decomposition and composition of deep CNNs, as discussed in \cite{guEtAl2022, guEtAl2023}, offer insights into accelerating training through subnetwork transfer learning. In \cite{guntherEtAl2020}, the authors explore layer-parallel training of deep residual neural networks (ResNets), demonstrating their effectiveness in tasks like image classification. In \cite{kopanickovaEtAl2023}, the authors introduce preconditioning strategies for physics-informed neural networks (PINNs), employing a nonlinear layer-wise preconditioner for the LFBGS optimizer.
Finally, in \cite{klawonnEtAl2023}, the authors present a novel CNN-DNN architecture that is inspired by DD methods, specifically made for model-parallel training. This architecture is designed to handle various image recognition tasks effectively. It incorporates local CNN models or subnetworks that process either overlapping or nonoverlapping segments of the input data, such as sub-images, demonstrating its adaptability and efficiency in different image recognition scenarios. These studies collectively underscore the significance of DD methods in training NNs across different domains.


In this work, we introduce a variant of the APTS method~\cite{grossKrause2009} to the field of machine learning.
Unlike conventional NN-training methods, the unique aspect of our proposed algorithm is its ability to:
\begin{itemize}
    \item[\tiny\ding{108}] guarantee global convergence in deterministic settings,
    \item[\tiny\ding{108}] eliminate the need for extensive hyper-parameter tuning, 
    \item[\tiny\ding{108}] and enable parallelization.
\end{itemize}
\noindent Thus, APTS allows for a first step to perfrom training on large supercomputers with efficient and parameter-free approaches.
We demonstrate the convergence properties of the proposed method numerically by utilizing benchmark problems from the field of classification.
Our numerical study includes a comparative analysis with the SGD, Adam, and LFBGS methods.

\section{Training Methods}
In this section, we provide a comprehensive explanation of the fundamental principles and techniques associated with the TR and APTS methods, focusing on their application to the training of NNs.

\subsection{Trust-Region Method}
TR methods \cite{nocedalWright1999, connGouldToint2000} are advanced globally convergent iterative algorithms utilized in non-linear, convex and non-convex optimization to find local minima of objective functions. These methods iteratively construct a simplified model of the objective function within a localized vicinity of the iterate estimate, referred to as the \textit{trust region}. 
The size of the trust region is adaptively adjusted at each iteration to ensure an accurate approximation of the true function behavior.

More precisely, in the TR framework at iteration $k$, the optimization problem~\eqref{eq:min_problem} is approximated by a  first- or second-order model $m_k$.
In order to obtain the search direction~$s^k$, the model $m_k$ is minimized as follows:
\begin{equation}\label{eq:tr_subproblem}
\min_{\| s \| \leq \Delta^k} m^k(s) = \min_{\| s \| \leq \Delta^k} \nabla \ell(\theta^k; \mathcal{D})^T (\theta^k+s) + \frac{1}{2} (\theta^k+s)^T H^k (\theta^k+s)\,,
\end{equation}
where $H^k \in \mathbb{R}^{n \times n}$ is a symmetric matrix, typically an approximation to the Hessian of $\ell$ at $\theta^k$, and $\Delta^k > 0$ is the TR radius. 
In case of a first-order model~$m_k$, the quadratic term in \eqref{eq:tr_subproblem} vanishes.
The quality of the step $s^k$ obtained by solving~\eqref{eq:tr_subproblem} is assessed using the ratio
\begin{equation}\label{eq:accuracy}
\rho^k = \frac{\ell(\theta^k) - \ell(\theta^k + s^k)}{m^k(0) - m^k(s^k)},
\end{equation}
where $\ell(\theta^k) - \ell(\theta^k + s^k)$ is the actual reduction and $m^k(0) - m^k(s^k)$ is the predicted reduction. 
The TR radius and parameters are then updated based on the value of $\rho^k$.
If $\rho_k \approx 1$, the model~$m_k$ provides a good approximation of the objective function, so the TR radius can be increased and the parameters can be updated. 
Otherwise, the TR radius is shrunk and the parameters remain unchanged. 
We refer the reader to \cite{connGouldToint2000} for a more detailed explanation and the algorithm.

In this work, we adopt the Limited-memory Symmetric Rank 1 (LSR1) method \cite{erwayMarcia2015, nocedalWright1999} for approximating the Hessian of~$\ell$.
The LSR1 approach systematically generates an approximation of the Hessian matrix by incorporating Rank 1 updates derived from gradient differences at each iteration. 
As the number of updates accumulates and reaches the predetermined limit, $m$, the most recent update replaces the oldest one, ensuring that the total number of updates remains constant at $m$. 
Unlike the LBFGS method, which ensures a positive definite approximated Hessian, the LSR1 method allows for an indefinite Hessian approximation, potentially providing more accurate curvature information.

In order to solve efficiently the subproblem \eqref{eq:tr_subproblem}, we employ the Orthonormal Basis SR1 (OBS) technique proposed in~\cite{brust2017}.
The OBS method leverages an orthonormal basis spanned by the LSR1 updates to transform the dimensions of the TR subproblem from $n \times n$ to a more manageable $m \times m$, where typically $m \ll n$. 

\subsection{APTS with the decomposition of network's parameters}
APTS~\cite{grossKrause2009} is a right-preconditioned TR strategy designed for solving non-convex optimization problems. 
In this work, we construct the nonlinear preconditioner by utilizing the layer-wise decomposition of NN's parameters~\cite{kopanickovaEtAl2023}. 
Thus, in order to define the subdomains, we introduce the projectors
\begin{equation}\label{eq:disjoint_ddomain}
    R_i:\mathbb{R}^n \longrightarrow \mathbb{R}^{n_i},\ i=1,...,N,\qquad \text{with}\qquad  R_iR_j^T=0,\ \forall i\neq j,
\end{equation}
where the second property implies that the subsets of parameters generated by $R_i$ are pair-wise disjoint. 
Here, $n$ denotes the dimensionality of the global parameter vector, and $n_i$ represents the dimensionality of the local parameter vector, associated with the $i$-th subdomain. 
We can now define the $i$-th network copy as $\mathcal{N}_i(\cdot;\theta_{i},\phi_{i}, \mathcal{D})$, where the subdomain parameters $\theta_{i}=R_i \theta$ are defined as the set of trainable parameters.
Moreover, with a slight abuse of notation, we define $\phi_{i}=\theta\setminus\theta_{i}$ as the remaining non-trainable parameters. 
Clearly the overall structure of $\mathcal{N}$ does not change, what changes is the trainable/non-trainable setting of its parameters.

\begin{algorithm}
\caption{Scheme of APTS/SAPTS in weight}\label{alg:APTS_W}
\DontPrintSemicolon
\KwIn{Objective function $\ell$, the dataset $\mathcal{D}$ or mini-batches $D_1,D_2, \ldots$, main NN $\mathcal{N}(\cdot;\theta^0, \mathcal{D})$, NN copies $\mathcal{N}_i(\cdot;\theta_{i,0}^0,\phi^0, \mathcal{D})$, maximum local TR iterations $\nu$, the value of the Boolean \texttt{fdl} variable, global TR radius $\Delta_G$ and the min/max global TR radii $\Delta_{\min}^G/\Delta_{\max}^G$, the decrease/increase factors $0<\alpha<1<\beta$ (local and global TR), the min/max radii $\Delta_{\min}^L/\Delta_{\max}^L$ of the local TR, the reduction/acceptance ratios $0<\eta_1<\eta_1<1$ (local and global TR).}
\texttt{epoch} $\gets 0$\;
\While{not converged}{ \label{algline:conv_1}
$k \gets 0$\;  
  \For{$D \in [\mathcal{D}]$ \texttt{(APTS)} or $D \in [D_1,D_2, \ldots]$ \texttt{(SAPTS)}}{\label{algline:conv_2}
    $\ell^k \gets \ell(\theta^{k};D)$\; \label{algline:old_loss}
    $g^k \gets - \frac{\nabla \ell(\theta^{k};D)}{\|\nabla \ell(\theta^{k};D)\|_\infty}\Delta_G$\; \label{algline:gcomp}
    \For{$i \in \{1, \ldots, N\}$ in parallel}{
      Set $(\theta_{i,0}^k,\phi^k) \gets (R_i \theta^k, \theta^k \setminus R_i \theta^k)$ and $\Delta_{L,i} \gets \Delta_G$\; \label{algline:sync_1}
      $\theta_{i,\text{new}} \gets$  $ \min_{\theta_i} \ell(\theta_i;D)$ using $\nu$ TR steps on
      $\mathcal{N}_i(\cdot;\theta_{i,0}^k,\phi^k, D)$ \; \label{algline:local_TR} 
    }
    Gather and sum up all local updates: $\tilde{\theta} \gets \sum_{i=1}^N R_i^T \theta^k_{i,\text{new}}$\; \label{algline:aggregation}
    Evaluate the preconditioning step: $s^k \gets \theta^{k}-\tilde{\theta}$\; \label{algline:step}    
    Update the weights of the main model $\mathcal{N}$ as $\theta^{k+\frac{1}{2}} \gets \theta^k + s^k$\; \label{algline:precond_step}
    $w \gets 0, \tilde{\Delta}\gets\Delta_G$\; \label{algline:global_model_weight_update}
    \While{$\ell(\theta^{k+\frac{1}{2}};D) > \ell^k$ and $w \neq 1$ and $\texttt{fdl}==\texttt{True}$}{
        \label{algline:loss_insurance_begin}
      $w \gets \min\lbrace w+0.2,1\rbrace$ and $\tilde{\Delta} \gets \max\lbrace\Delta_{\min},\alpha\tilde{\Delta}\rbrace$\;
      $\theta^{k+\frac{1}{2}} \gets \theta^k + \frac{wg^k+(1-w)s^k}{\|wg^k+(1-w)s^k\|_\infty}\tilde{\Delta}$\; \label{algline:loss_insurance_end}
    }
    $\theta^{k+1} \gets$  $ \argmin_{\theta} \ell(\theta;D)$ using one TR step on $\mathcal{N}(\cdot,\theta^{k+\frac{1}{2}},D)$\; \label{algline:global_pass}
    $k \gets k + 1$\; 
  }
  $\texttt{epoch} \gets \texttt{epoch} +1$\;
}
\end{algorithm}

APTS algorithm can be executed using the entire dataset or on minibatches, with the latter scenario being referred to as Stochastic APTS (SAPTS). 
Although the convergence proof of the deterministic version does not transfer directly to the stochastic version, the numerical results provided in Section \ref{sec:num_results} are promising. 
Algorithm \ref{alg:APTS_W} provides a summary of the proposed APTS training algorithm. 
In steps \ref{algline:conv_1} and \ref{algline:conv_2}, we establish the convergence criteria and select the training approach, which can be either using the full dataset (APTS) or employing a mini-batch strategy (SAPTS). 
In steps \ref{algline:old_loss}-\ref{algline:gcomp}, we compute and store the loss and gradient with the current training dataset. The vector $g$ is normalized in step \ref{algline:gcomp} so that it has a similar length to $s$ (see step \ref{algline:step}), such that their contribution is similar in step \ref{algline:loss_insurance_end}. 
In step \ref{algline:sync_1}, we synchronize all the local TR radii $\Delta_{L,i}$ with the global one and all the network copies, $\mathcal{N}_i$, with the current global parameters from $\mathcal{N}$. This synchronization preserves the state of trainable and non-trainable parameters that were previously established.
Following this, in step \ref{algline:local_TR}, we initiate a series of TR iterations for each local network copy, which we call \textit{local TR iterations}. 
This process iterates until the cumulative step size falls below the predefined global TR radius, $\Delta_G$, signified by the condition $s^k_{i,j}=\|\theta_{i,0}^k-\theta_{i,j}^k\|_\infty\leq\Delta_G$, where $i$ denotes the submodel and $0$ and $j$ are the local TR iterations. 
The iteration halts when either $s^k_{i,j}$ equals $\Delta_G$ or when the iteration count reaches $\nu$.

Following the local TR iterations, step \ref{algline:aggregation} involves aggregating all the updated weights from the local networks and summing them up. Here, we emphasize that the update $s^k$ made in step \ref{algline:step} is bounded by the global TR radius $\Delta_G$, i.e.,
\begin{equation*}
    \| \theta^k - \tilde{\theta}\|_\infty
    =\| \sum_{i=1}^N R_i^T\theta^k_{i,0} - \sum_{i=1}^N R_i^T \theta^k_{i,\text{new}}\|_\infty
    =\| \sum_{i=1}^N R_i^T(\theta^k_{i,0}-\theta^k_{i,\text{new}})\|_\infty \leq \Delta_G.
\end{equation*}
In addition, the operator $R_i^T$ serves the specific purpose of padding the vector $\theta^k_{i,\text{new}}$ with zeros, transforming it to the same dimensionality as $\theta^k$. This ensures a non-overlapping aggregation during the summation process of $R_i^T\theta^k_{i,\text{new}}$ across $i=1,...,N$, a condition which follows from the property of $R_i$ delineated in equation \eqref{eq:disjoint_ddomain}. 

In step \ref{algline:precond_step}, we evaluate the preconditioning step, i.e., the search direction obtained during the preconditioning iteration. 
This direction is then used to update the weights of the global model in step \ref{algline:global_model_weight_update}. 
In steps \ref{algline:loss_insurance_begin}-\ref{algline:loss_insurance_end}, we make sure that the current loss $\ell(\theta^{k+\frac{1}{2}}; D)$ on the global model is smaller than the old loss $\ell^k$.
If this is not a case, then we modify $s^k$ by averaging it with the gradient $g^k$. 
Moreover, we reduce the global TR radius $\Delta_G$. 
This process is different from the one in~\cite{grossKrause2009} and it is inspired by the Dogleg method \cite{nocedalWright1999} which is repeated until a step is accepted or until the step is the gradient, i.e. when $w=1$.\\
In the final phase, a single TR iteration is executed on the global neural network \(\mathcal{N}\) in step \ref{algline:global_pass}.

\section{Numerical Results}\label{sec:num_results}
In this section, we compare the performance of APTS against well-established training algorithms such as SGD, LBFGS and Adam. All tests, accessible through the repository in \cite{cruzas_ML_APTS}, were conducted using Python 3.9, with the implementation leveraging PyTorch version 2.0.1+cu117. For SGD, LBFGS, and Adam, we utilized the standard optimizers provided by PyTorch.

We tuned the hyper-parameters of SGD, LBFGS, and Adam to ensure a fair comparison, recording the average best results in 10 trials across various parameter configurations. Specifically:
\begin{itemize}
    \item[\tiny\ding{108}] For SGD, we evaluated 63 combinations, varying the learning rate in $[10^{-4},1]$ and the momentum in $[0.1,0.9]$.
    \item[\tiny\ding{108}] For LBFGS, we considered 45 combinations, adjusting both the learning rate within the interval $[10^{-4},1]$ and the history size (the amount of stored updates of the Hessian) within $[3,50]$. Note that we set the \texttt{line\_search\_fn} parameter to \texttt{strong\_wolfe}. 
    \item[\tiny\ding{108}] For Adam, we examined 13 different learning rates within the interval $[10^{-5},10^{-1}]$.
\end{itemize}
For these three algorithms, all other input parameters were kept as their standard PyTorch default values. Moreover, the selected optimal parameters are shown in the legend of the  figures. In contrast, APTS and TR were consistently run with the same set of input parameters across all tests. Detailed descriptions of these settings can be found in the next section. 

Afterwards, we provide a detailed analysis of the obtained results, highlighting the characteristics of APTS in comparison to SGD, LBFGS, Adam and TR. Note that in the absence of the mini-batch strategy, the SGD optimizer simplifies to the traditional Gradient Descent (GD) algorithm. We remark that one epoch of APTS consists in the preconditioning step and the global step in lines \ref{algline:conv_2}--\ref{algline:global_pass} of Algorithm~\ref{alg:APTS_W}.

Finally, some of the following plots will depict the average test accuracy and training loss (represented by a solid line) alongside the variance observed in the worst and best runs (indicated by dashed lines) across 10 independent trials, each initialized with different random NN parameters.

\subsection{APTS Implementation Details and Parameter Configuration}\label{subsec:implementation_details}

The PyTorch programming environment imposes certain constraints on the level of control one has over parameter training. Parameters are configured in vectors, with entire vectors designated as either trainable or non-trainable. This means that, for instance, a fully connected layer, which consists of weights and biases, is stored as a structure composed by two distinct vectors of learnable parameters, and they can be configured to be trained independently. Therefore, the projectors $R_i$ in equation~\eqref{eq:disjoint_ddomain} are built in a way to maximize the efficiency in training, i.e., there will not be two weights belonging to the same PyTorch vectors with different trainable/non-trainable settings. Furthermore, the projectors are randomly generated, i.e., if we have a total of $M$ vectors of parameters and we want to split them across $N$ models, each model will have approximately round$(\frac{M}{N})$ randomly selected trainable vectors. 

Key configurations in our implementation of the APTS algorithm are as follows:

\begin{itemize}
    \item[\tiny\ding{108}] In step 8 of APTS, i.e. Algorithm \ref{alg:APTS_W}, the amount of subdomains (submodels $\mathcal{N}_i$) is specified in the legend of the figures in Section \ref{sec:num_results} and the trainable weight vectors are randomly selected.
    \item[\tiny\ding{108}] In step 10, concerning the \textit{local} TR, we configured:
    \begin{itemize}
        \item[\textbf{-}] A maximum amount of iterations to $\nu=5$;
        \item[\textbf{-}] A minimum radius of $\Delta^L_{\min}=0$, an initial radius equal to the current global TR radius $\Delta^L=\Delta^G$, and a dynamically adjusted maximum radius $\Delta^L_{\max}=\Delta^G-\|\theta_{i,0}^k-\theta_{i,j}^k\|$ at the $j$-th iteration to ensure that the cumulative step remains bounded by the global TR radius;
        \item[\textbf{-}] Whenever a quadratic model \eqref{eq:tr_subproblem} is considered, history size for the SR1 Hessian is set to $5$. Moreover, at the start of the preconditioning iteration, the local Hessian history is synchronized with the global one;
        \item[\textbf{-}] The reduction and acceptance ratios $0<\eta_1<\eta_2<1$ are the standard values for TR, i.e., $\eta_1=0.25$ and $\eta_2=0.75$;
        \item[\textbf{-}] The decrease and increase factors $\alpha$ and $\beta$ are set to $0.5$ and $2$, respectively.
    \end{itemize}
    \item[\tiny\ding{108}] In step 18, which involves the \textit{global} TR:
    \begin{itemize}
        \item[\textbf{-}] Minimum, maximum, and initial radii are set to $\Delta^G_{\min}=10^{-4}$, $\Delta^G_{\max}=10^{-1}$, and $\Delta^G=10^{-2}$, respectively.
        \item[\textbf{-}] The values of $\alpha,\beta,\eta_1,\eta_2$ match the ones of the local TR;
        \item[\textbf{-}] The reduction and acceptance ratios match those of the local TR;
        \item[\textbf{-}] Whenever a quadratic model \eqref{eq:tr_subproblem} is considered, history size for the SR1 approximated Hessian is set to $5$.
    \end{itemize}
\end{itemize}
Note that both the global and local TR use the infinity norm to compute the length of the steps. This does not conflict with 2-norm used by the OBS method, as the 2-norm is bounded by the infinity-norm. 

In the numerical results Section \ref{sec:num_results}, the TR method as standalone solver, denoted by TR in the figures, has the same settings as the global TR method in APTS.

\subsection{MNIST}
We conducted some preliminary experiments using the MNIST dataset \cite{lecun1998}. This dataset comprises 60'000 training images and 10'000 test images of handwritten digits, each of which is a $28\times28$ pixel grayscale image. 
For this experiment we considered a Fully Connected Neural Network (FCNN) trained on the entire dataset, i.e., without the use of a mini-batch strategy.

The FCNN used in our experiments consists of three fully connected layers (6 PyTorch vectors), resulting in a total of 26'506 parameters. To investigate the performance of APTS, we configured it with 2, 4, and 6 subdomains and with an enabled/disabled forced-decreasing-loss parameter (\texttt{fdl}), while other parameters are reported in Section \ref{subsec:implementation_details}.

\begin{figure}[ht]
\centering
\begin{subfigure}{0.49\textwidth}
  \centering
  \includegraphics[width=\linewidth]{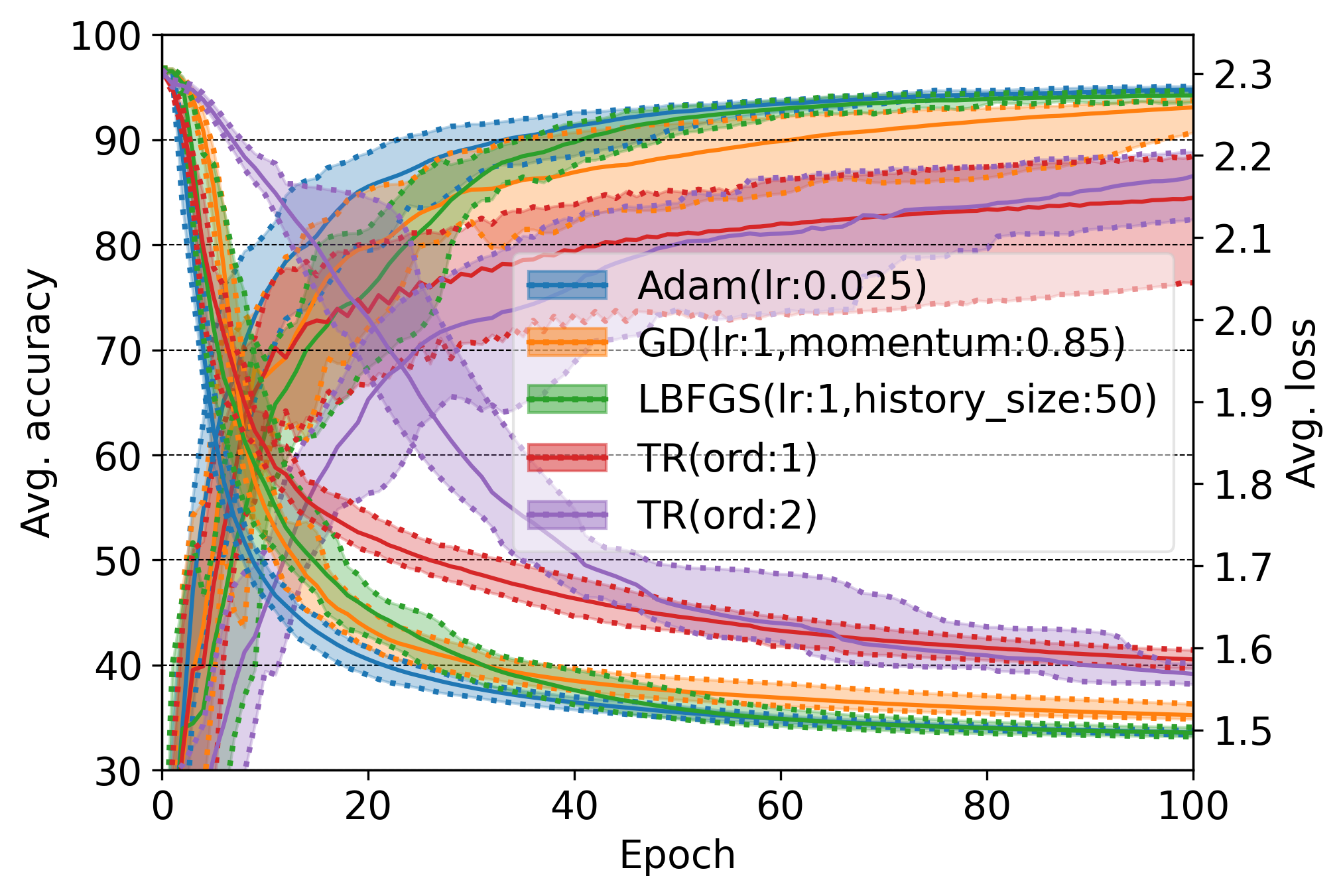}
  \caption{Standard baseline optimizers}
  \label{fig:sub1}
\end{subfigure}
\hfill
\begin{subfigure}{0.49\textwidth}
  \centering
  \includegraphics[width=\linewidth]{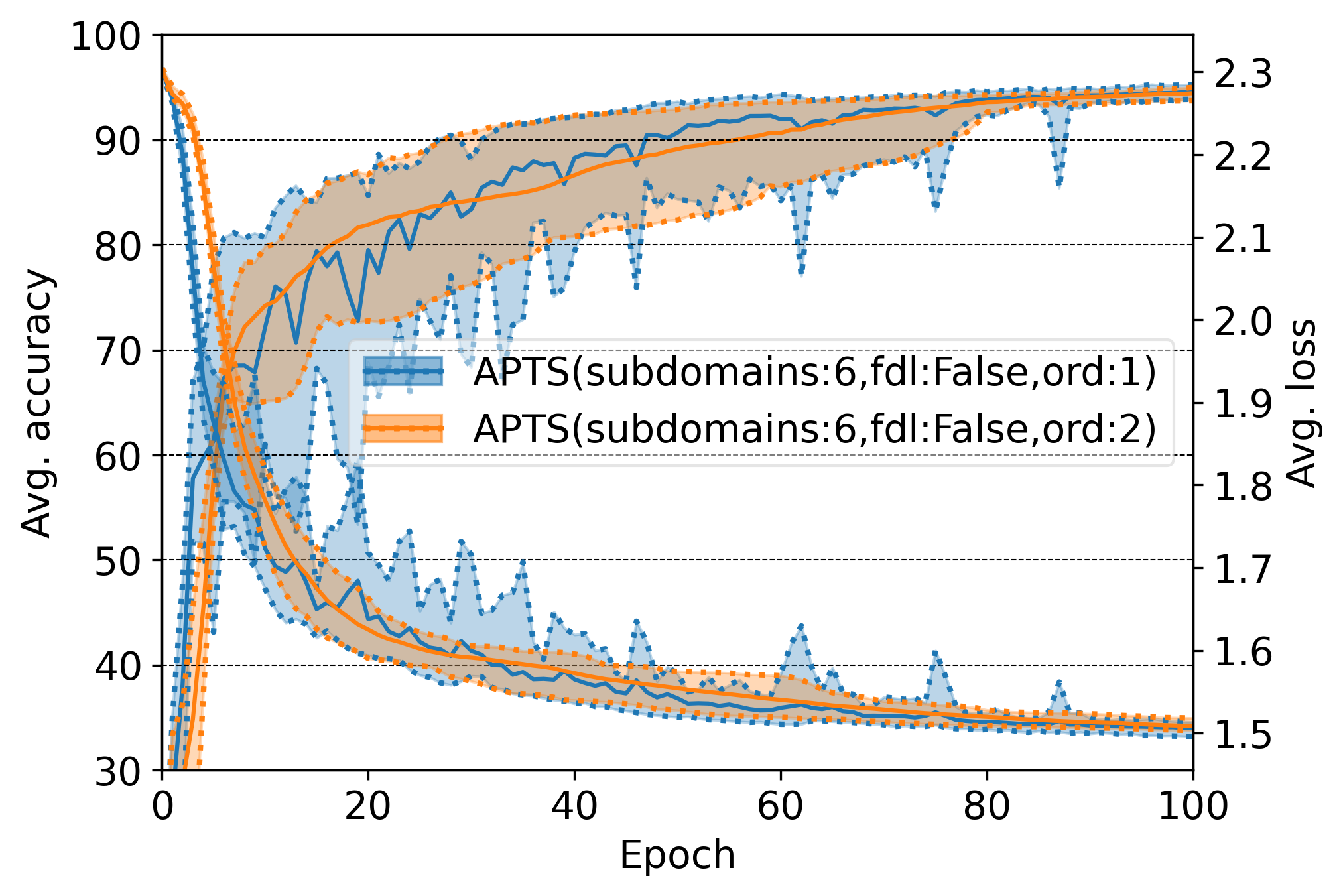}
  \caption{APTS without \texttt{fdl}}
  \label{fig:sub2}
\end{subfigure}

\begin{subfigure}{0.49\textwidth}
  \centering
  \includegraphics[width=\linewidth]{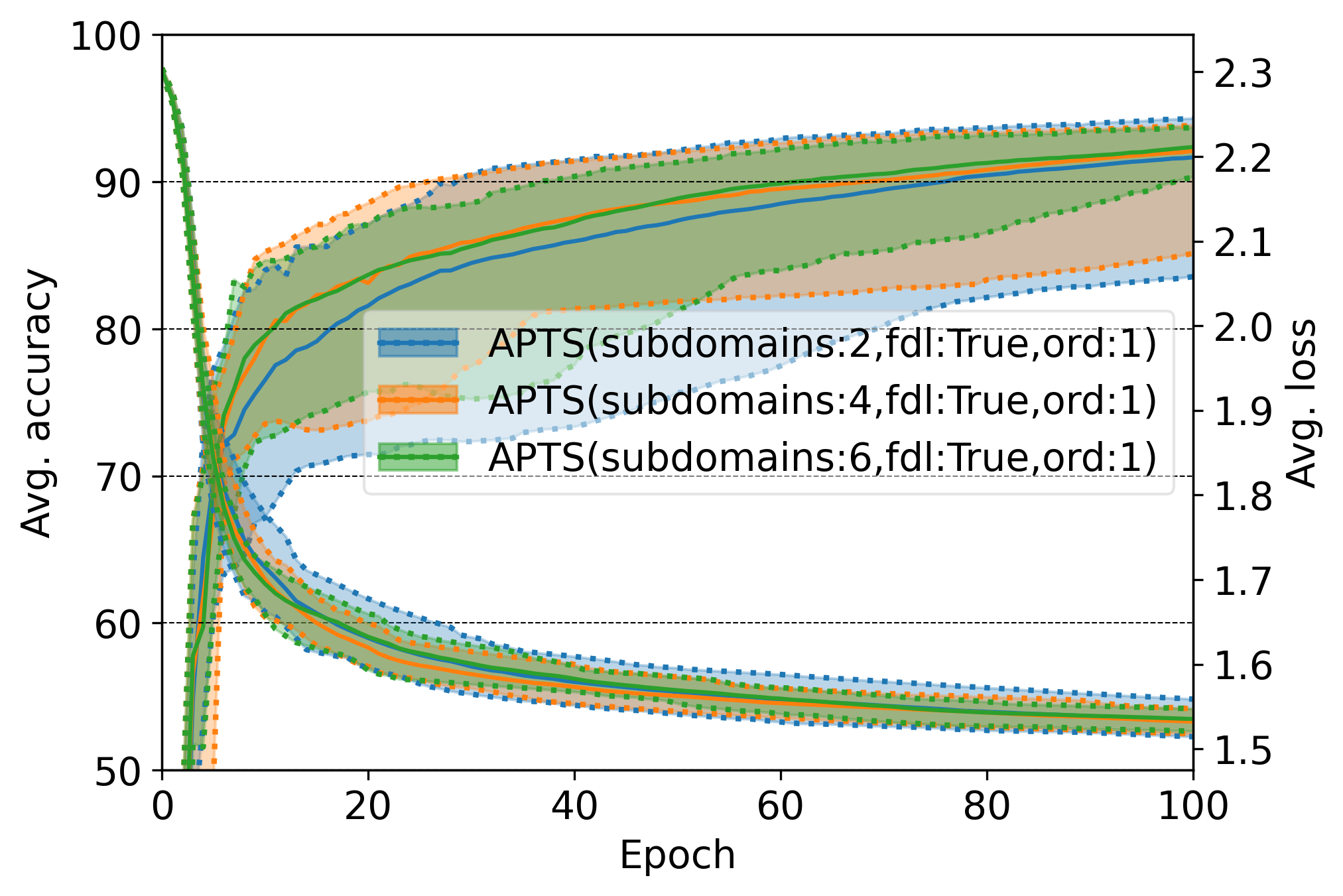}
  \caption{APTS with first-order TR subproblem}
  \label{fig:sub3} 
\end{subfigure}
\hfill
\begin{subfigure}{0.49\textwidth}
  \centering
  \includegraphics[width=\linewidth]{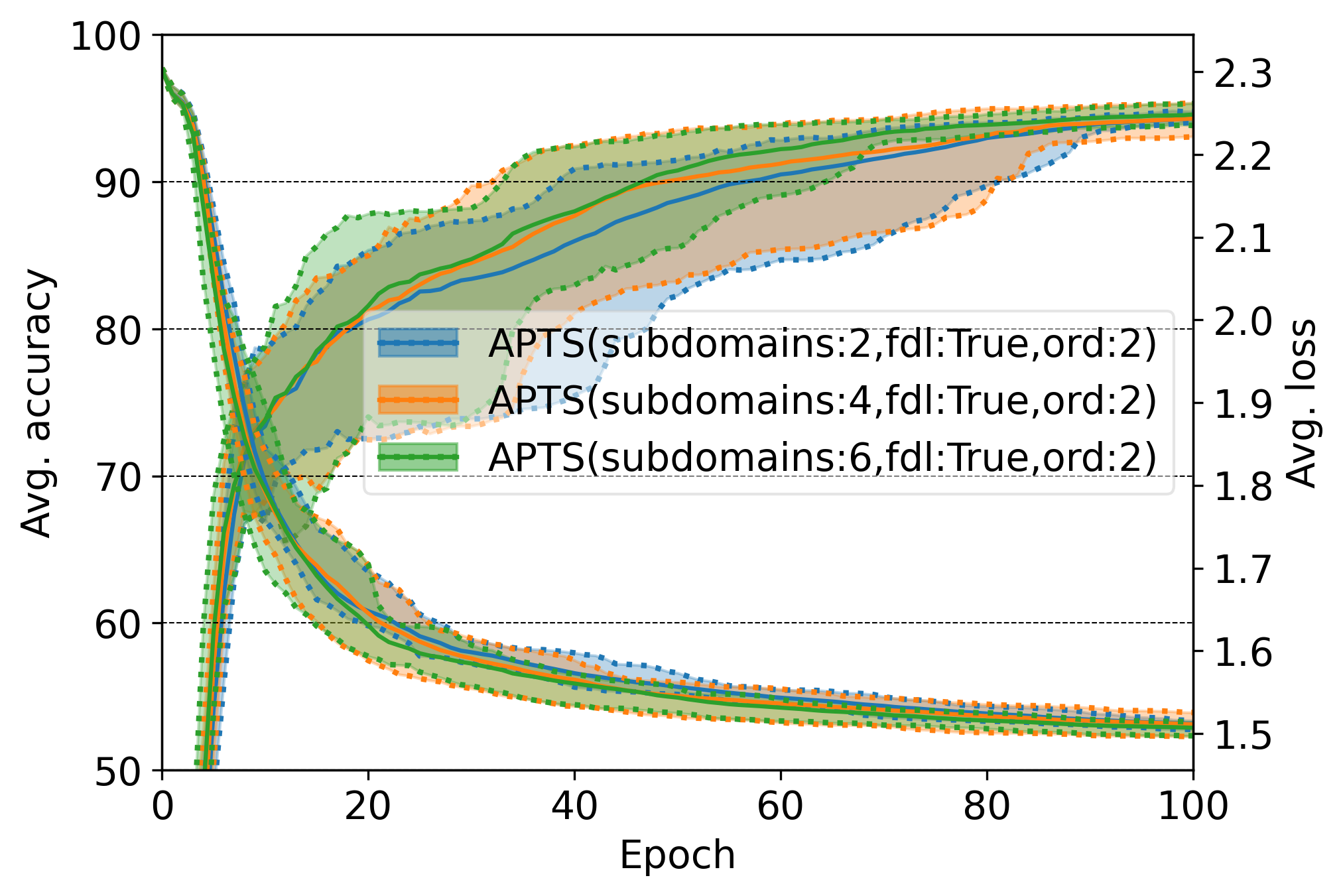}
  \caption{APTS with second-order TR subproblem}
  \label{fig:sub4} 
\end{subfigure}
\caption{Comparison of APTS variants and baseline optimizers on FCNN trained on MNIST.}
\label{fig:fcnn_results}
\end{figure} 

Observing Figure \ref{fig:sub1}, we note that Adam, LBFGS, and GD achieve comparable levels of accuracy. However, GD exhibits a marginally higher number of epochs to converge and a greater variance in its results. The TR variants, utilizing both first- and second-order information in equation \eqref{eq:tr_subproblem} (denoted by $\texttt{ord}=1$ and $\texttt{ord}=2$, respectively), demonstrate a slower convergence compared to the benchmark algorithms. TR($\texttt{ord}=2$) overtakes its first-order counterpart TR($\texttt{ord}=1$) after 70 epochs. This can be attributed to the enhanced reliability of the approximated Hessian close to a minimizer of the training loss, which has a flatter landscape in comparison to the initial training phase. Moreover, the introduction of second-order information tends to stabilize the variation among TR iterations compared to the first-order approach.

In Figure \ref{fig:sub2}, the APTS variants employing first- and second-order derivatives with \texttt{fdl} turned off are presented. The use of second-order derivatives seems to facilitate a smoother progression of the loss function over the training epochs. Conversely, reliance on first-order derivatives results in more fluctuating loss values. Nonetheless, both variants converge to comparable levels of loss and accuracy, with a limited variance suggesting a consistent convergence pattern across different runs.

Figures \ref{fig:sub3} and \ref{fig:sub4} showcase the APTS variants with \texttt{fdl} enabled. Here, the application of second-order information leads to reduced variance and marginally superior accuracy relative to the first-order variant. Additionally, an increase in the number of subdomains correlates with improved performance. APTS with 6 subdomains not only slightly outperforms the 2 and 4 subdomain configurations in terms of accuracy, but also shows a reduced variance. This enhanced performance likely results from the added subdomains enabling more thorough exploration of the loss landscape through independent parameter adjustments, in contrast to the constrained exploration provided by a lower count of subdomains. This concept aligns with prior studies advocating for distinct weight and bias training strategies to improve training dynamics \cite{ainsworthShin2021}.

A comparison of Figures \ref{fig:sub3} and \ref{fig:sub2} reveals that APTS with six subdomains, with \texttt{fdl} disabled and first-order information, attains a marginally superior accuracy and lower loss than its counterpart with \texttt{fdl} enabled. This highlights the importance of unrestricted parameter exploration in improving training efficiency, similar to the exploration afforded by a higher subdomain count.

\subsection{CIFAR-10}
Following our MNIST investigation, we evaluated the solution strategies using the CIFAR-10 dataset \cite{krizhevsky2009learning}, containing 50'000 training and 10'000 test $32\times32$ color images across $10$ classes. For these experiments, we employed the ResNet18 architecture from the PyTorch library with 11'689'512 parameters stored in 62 PyTorch vectors. Training was conducted using a mini-batch strategy, with $20$ overlapping mini-batches, each containing 2'975 samples, which were created with a $1\%$ overlap between the mini-batches in order to reduce variance.

\begin{figure}[ht]
\centering
\begin{subfigure}{0.49\textwidth}
  \centering
  \includegraphics[width=\linewidth]{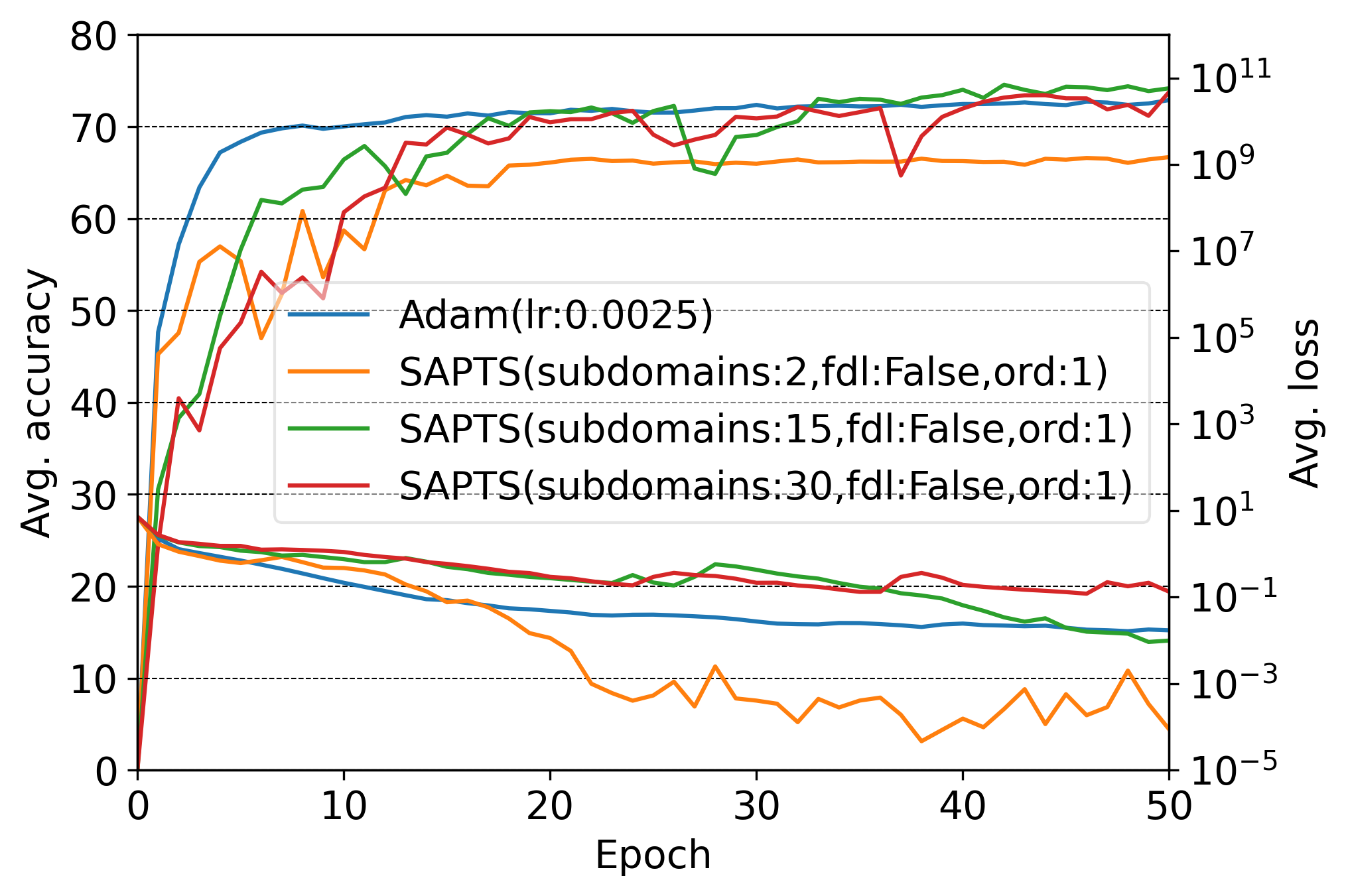}
  \caption{Adam and APTS with $2,15,30$ subdomains}
  \label{fig2:sub1}
\end{subfigure}
\hfill
\begin{subfigure}{0.49\textwidth}
  \centering
  \includegraphics[width=\linewidth]{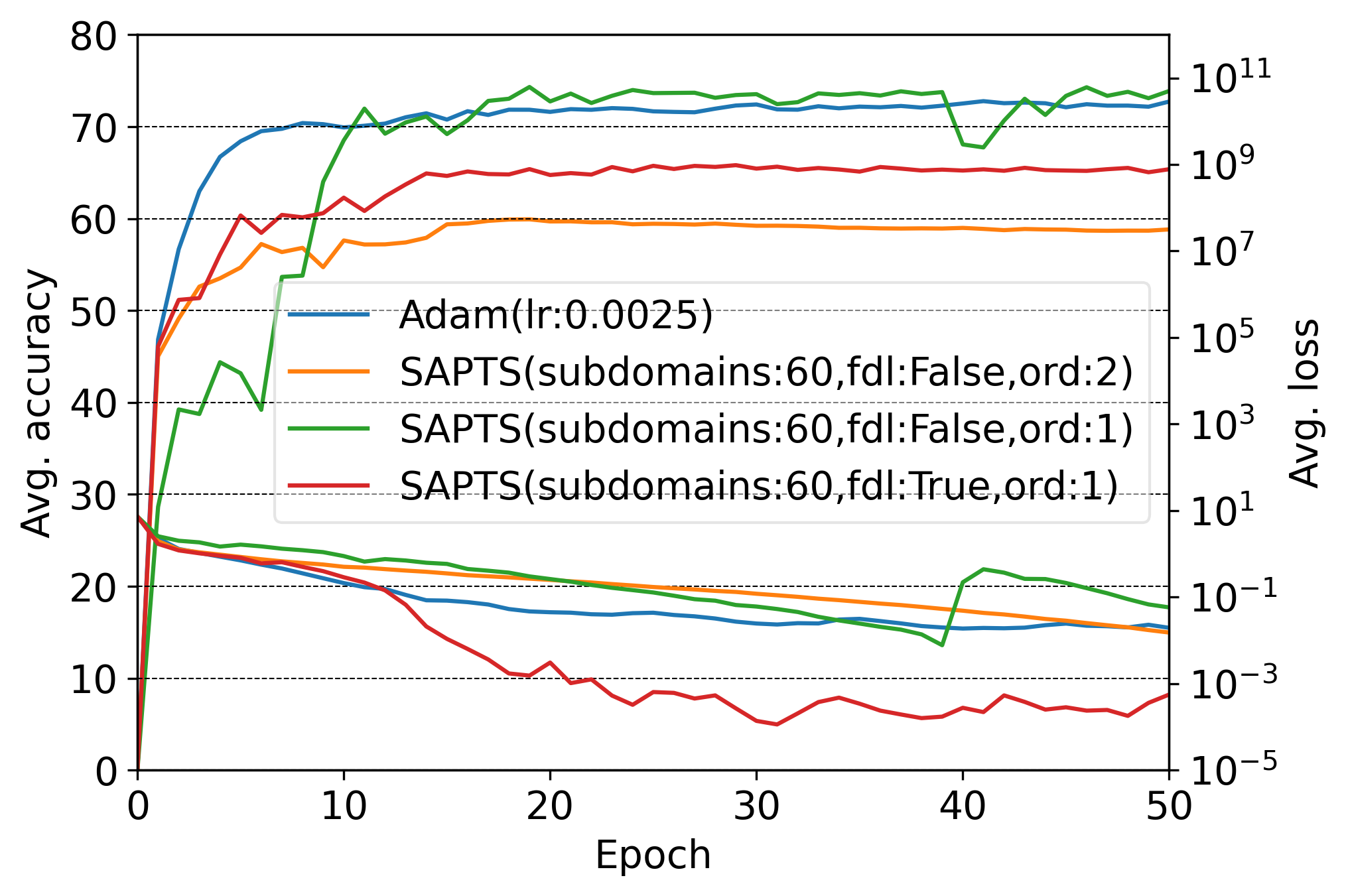}
  \caption{Adam and APTS variants with $60$ subdomains}
  \label{fig2:sub2}
\end{subfigure}
\caption{Comparison of APTS variants and baseline optimizers on CNN trained on CIFAR-10.}
\label{fig2}
\end{figure}

In Figure \ref{fig2}, we present a comparative analysis between stochastic Adam and various configurations of SAPTS. Specifically, Figure \ref{fig2:sub1} features SAPTS with 2, 15, and 30 subdomains, a first-order TR approach for the subproblems, and with the feature \texttt{fdl} deactivated. In Figure \ref{fig2:sub2}, the number of subdomains is set to 60 while we vary parameters such as \texttt{fdl} and \texttt{order}. We do not provide results for SGD here as Adam achieved a higher test accuracy, and our aim is to provide a comparison with APTS.

Our primary focus was on configurations with $\texttt{fdl}=\texttt{False}$ and $\texttt{ord}=1$. The rationale behind this is that activating \texttt{fdl} (as shown in Figure \ref{fig2:sub2}) leads to a deterioration in the generalization abilities of SAPTS, albeit with a higher reduction in loss compared to other SAPTS variants. Regarding $\texttt{ord}=2$, SAPTS exhibited poor generalization and did not significantly reduce the loss. However, it maintained a smoother loss curve compared to other tested SAPTS variants as in the previous non-stochastic cases.

When analyzing SAPTS with $\texttt{fdl}=\texttt{False}$ and $\texttt{ord}=1$ across different numbers of subdomains, we observe a correlation similar to that in non-stochastic tests: SAPTS's performance improves with an increased number of subdomains. Additionally, when comparing SAPTS with 60 subdomains against Adam, we find that SAPTS demonstrates marginally better generalization over a longer period and requires only a few more epochs to match the performance level of Adam.

\section{Conclusion}
In this work, we employed APTS method, which utilized the decomposition of the NN's parameters, to efficiently train neural networks.  
Our numerical study reveals that the proposed APTS method 
attains 
comparable or superior generalization capabilities to traditional optimizers like SGD, LBFGS, and Adam. 
The strength of the developed APTS method lies in absent hyper-parameter tuning and in its parallelization capabilities. 
Moreover, the APTS achieves 
faster convergence with a growing number of subdomains, showing encouraging algorithmic scalability. 


\section*{Acknowledgement}
This work was supported by the Swiss Platform for Advanced Scientific Computing (PASC) project ExaTrain (funding periods 2017-2021 and 2021-2024) and by the Swiss National Science Foundation through the projects ``ML$^2$ - Multilevel and Domain Decomposition Methods for Machine Learning'' (197041) and ``Multilevel training of DeepONets - multiscale and multiphysics applications'' (206745).

\printbibliography

\end{document}